\documentclass[12pt]{article}
\usepackage{hyperref}
\usepackage{graphicx}
\usepackage{amsfonts}
\usepackage{amsmath}
\usepackage{tabularx, booktabs}
\usepackage{siunitx}
\usepackage{array}
\usepackage{color,soul}
\usepackage{xcolor}
\usepackage[T1]{fontenc}
\usepackage{mathtools}

\DeclarePairedDelimiter\floor{\lfloor}{\rfloor}
\newcolumntype{M}[1]{>{\centering\arraybackslash}m{#1}}
\begin{document}
\title{Definitive Proof of Goldbach's Conjecture}
\author{Kenneth A. Watanabe, PhD}
\maketitle

\section{Abstract}

The Goldbach conjecture states that every even integer greater than 2 can be expressed as the sum of two prime numbers. This conjecture was first proposed by German mathematician Christian Goldbach in 1742 and, despite being obviously true, has remained unproven. In this paper, it is shown that the set of all even integers $n$ that are not divisible by a prime number less than $\sqrt{n}$ has the relatively fewest number of prime pairs. An equation was derived that approximates the number of prime pairs for these values of $n$. It was then proven that this equation never goes to zero for any $n$, and as $n$ increases, the number of prime pairs also increases, thus validating Goldbach's conjecture. Error analysis was performed to show that the difference between this approximation and the actual number of prime pairs is small enough so that for all $n>622$, the number of prime pairs of $n$ is greater than 1, thus proving Goldbach's conjecture.

\section{Introduction}

On June 7th, 1742, the German mathematician Christian Goldbach \cite{Goldbach} wrote a letter to Leonhard Euler in which he proposed that every even integer greater than 2 can be expressed as the sum of two primes. Since the time it was proposed, the conjecture was widely believed to be true - Euler himself replied to Goldbach: "That ... every even integer is the sum of two primes, I regard as a completely certain theorem, although I cannot prove it." \cite{Ingham}.  However, the conjecture has remained unproven for over 250 years. In March 2000, the publisher Faber and Faber Co. offered a \$1 million prize to anyone who can prove Goldbach's conjecture. Unfortunately, the offer expired in 2002 and the prize went unclaimed.

\begin{figure}
  \includegraphics[width=\linewidth]{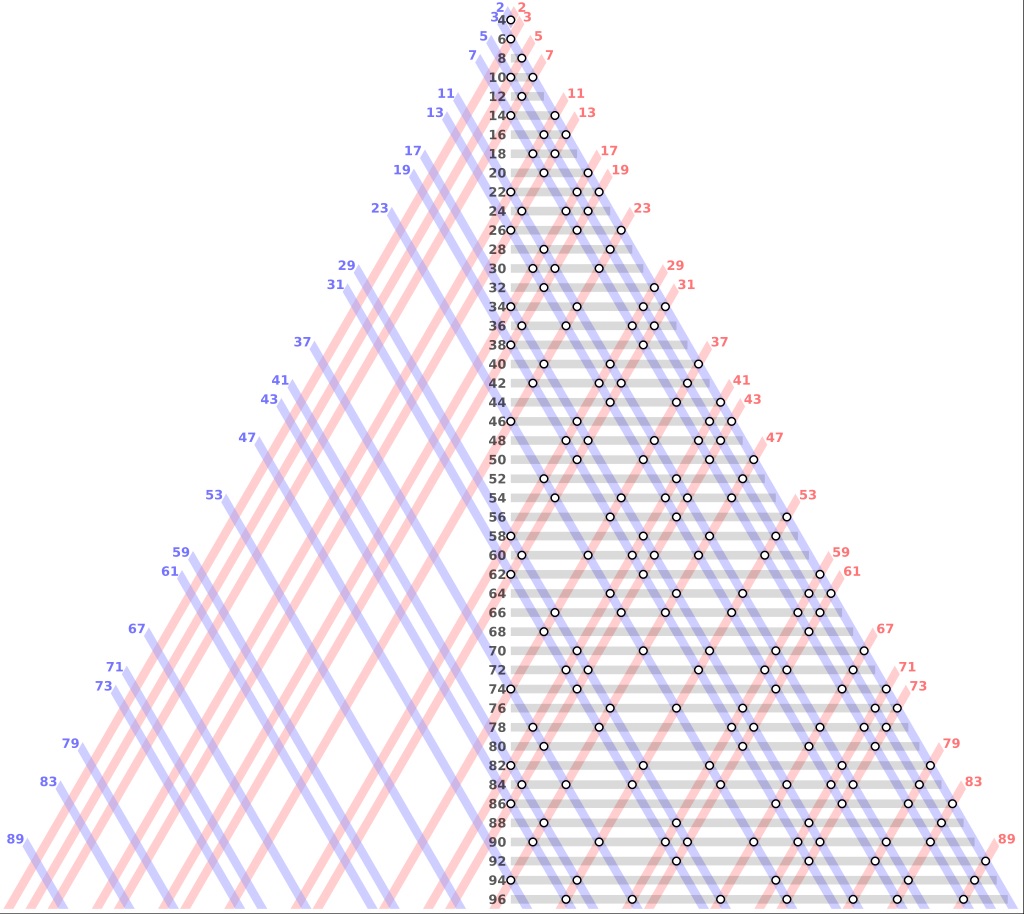}
  \caption{Goldbach partitions of even numbers from 4 to 96. The open circles on the intersections of the purple and pink lines represent the prime pairs that total the $n$ value on the vertical column. \cite{Cunningham}}
  \label{fig:Goldbach_Partitions}
\end{figure}

A graphic representation of the Goldbach partitions by Cunningham and Ringland \cite{Cunningham} (Figure \ref{fig:Goldbach_Partitions}) shows all the prime pairs for even $n$ from $n = 4$ to $n = 96$. In general, the number of prime pairs increases as $n$ increases, but there are many exceptions. Notably, where $n = 2p$ (where $p$ is prime), the number of prime pairs is fewer than when $n=2\times3\times p$ or $n=2\times3\times5\times p$. For example, 94 has 5 prime pairs, but 96 has 7 prime pairs and 90 has 9 prime pairs. Same with 62 versus 66 or 60. In fact, it turns out that the values where $n$ has the fewest prime pairs occur when $n$ is not divisible by a prime number less than $\sqrt{n}$. This occurs when $n=2^i$ or $n=2^ip$ where $i$ is an integer greater than 0, and $p$ is a prime number greater than $\sqrt{n}$. This will be further discussed later in this paper.

\begin{figure}
  \includegraphics[width=\linewidth]{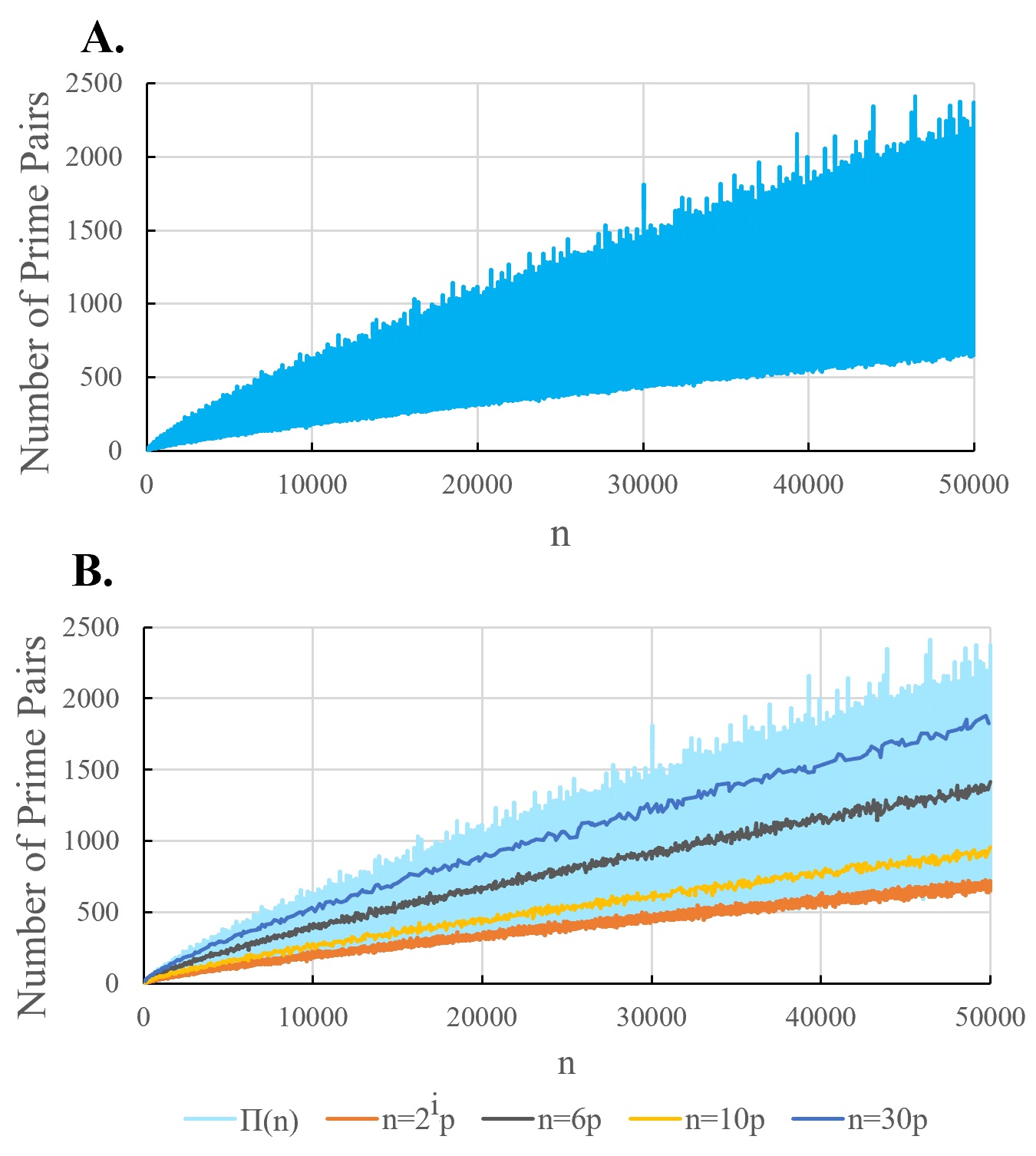}
  \caption{Graph of (A.) the number of prime pairs for even integers $n$ less than 50,000 (light blue).
  The lower bound on the number of prime pairs (B.) occurs when $n$ is not divisible by primes less than $\sqrt{n}$ (orange line) (i.e. $n=2^i$ or $n=2^ip$ where $p > \sqrt{n}$). For values of $n = 6p$ (black line), $n = 10p$ (yellow line) and $n = 30p$ (dark blue line), the relative number of prime pairs is larger. }
  \label{fig:Graphs}
\end{figure}

A graph of the prime pairs of even integers $n$ for values of $n$ up to 50,000 was made as shown in Figure \ref{fig:Graphs}. The curve appears to be increase indefinitely. Also graphed are the number of prime pairs for $n$ not divisible by a prime less than $\sqrt{n}$ (or $n=2^i$ or $n=2^ip$ where $p > \sqrt{n}$), $n=2\times3p$ (or $n = 6p$), $n = 2\times5p$ (or $n = 10p$) and $n = 2\times3\times5p$ (or $n = 30p$) where $p$ is a prime number (Figure \ref{fig:Graphs}B). The number of prime pairs is significantly fewer when $n$ is not divisible by a prime less than $\sqrt{n}$ (orange line) compared to the other values of $n$. The more prime factors of $n$, the higher the number of prime pairs as can be seen by the other curves $n=2\times3p$ (black line), $n=2\times5p$ (yellow line) and $n=2\times3\times5p$ (dark blue line). Note that all the curves have a positive slope and the gap between the curves increases with increasing $n$.

Goldbach's conjecture was computationally proven to be valid for all values of $n$ up to $n = 4 \times 10^{18}$ by Silva et al \cite{Silva} and the curves in their research look similar to those in Figure \ref{fig:Graphs}. This of course does not mean that there is no possibility that somewhere even further out that the curve will turn around and go down to zero or that maybe there is a single value for $n$ somewhere that cannot be expressed as the sum of two primes.
The following proof confirms that Goldbach's conjecture is true for all even integers $n$.

\section{Functions}
Before we get into the proof, let us define some functions that are necessary for the proof.

Let the function $l(n)$ represent the largest prime number less than $n$. For example, $l(10) = 7$, $l(20) = 19$ and $l(19.5) = 17$.

Let $P(n)$ represent the set of pairs $(x,y)$ such that $x+y = n$ and both $x$ and $y$ are odd integers greater than 2. The values of $x$ or $y$ need not be prime. Let $|P(n)|$ be the number pairs in $P(n)$, thus $|P(n)| = (n/2)-2$.

Let the function $K(n)$ represent the number of pairs $(x,y)$ such that $x+y = n$ and either $x$ or $y$ are odd composite integers greater 2.

Let the function $K^*(n)$ represent an approximation of $K(n)$ as $n \lim\to\infty$.

Let the function $W(n)$ represent the fraction of pairs $(x,y)$ where $x+y = n$ and both the x and y coordinates are prime numbers.

Let the function $\Pi(n)$ represent the actual number of pairs $(x,y)$ such that $x+y = n$ and both $x$ and $y$ are prime numbers greater than 2.

Let the function $\Pi^*(n)$ represent the approximation of $\Pi(n)$ as defined by the equation $\Pi^*(n) = |P(n)|W(n)$.

Let the function $\pi(n)$ represent the prime counting function. $\pi(n)$ is the number of prime numbers less than or equal to $n$.

The term ``number of pairs of n'' means the number of pairs (x,y) such that x+y=n and x and y are odd integers and n is an even integer.
The term ``a pair is divisible by p'' means that either the x or y coordinate of a pair is evenly divisible by p.  
The term ``a pair is not divisible by p'' means that both the x and y coordinates of a pair are not evenly divisible by p.

\section{Methodology for Proving Goldbach's Conjecture}

To prove Goldbach's conjecture, we have to show that for any given even integer $n$ there is at least one pair $(p_{x}, p_{y})$ such that both $p_{x}$, and $p_{y}$ are prime and $p_{x} + p_{y} = n$. Since all primes except 2 are odd, we will only consider odd integers for x and y, and we will exclude any pair containing 1 since 1 is not considered a prime number.
The first step is to find the set of all pairs of odd integers $(x,y)$ such that $x + y = n$ and $x$ and $y$ are both greater than 2. Then we must calculate the number of pairs $(u,v)$ such that either $u$ or $v$ is evenly divisible by a number > 2 (i.e. not prime), and that $u + v = n$ and u and v are odd integers.  If we can show that the number of pairs $(x,y)$, is greater than the number of pairs $(u,v)$, then we know that there exists at least one pair $(x',y')$ such that $x' + y' = n$ and $x'$ and $y'$ are prime numbers.\\
For every even integer $n$, there are exactly $n/2$ pairs of odd numbers $(x,y)$ such that $x + y= n$. Excluding pairs $(1,n-1)$ and $(n-1,1)$, there are exactly $|P(n)| = (n/2)-2$ pairs of odd numbers $(x,y)$ such that $x + y= n$ and $x$ and $y$ are both greater than 1.\\

The following table shows some examples of the values of $n$ and the number of pairs.

\begin{table}[h!]
\resizebox{\textwidth}{!}{%
\begin{tabular}{ |c|M{1.5cm}|p{1cm}|l| } 
\hline
 n&|P(n)|&$\Pi(n)$&P(n), Set of pairs (x,y) where x+y=n and x and y are odd integers > 2\\
\hline
6&1&1&(3,3)\\
8&2&2&(3,5), (5,3)\\
10&3&3&(3,7), (5,5), (7,3)\\
12&4&2&(3,9), (5,7), (7,5), (9,3)\\
14&5&3&(3,11), (5,9), (7,7), (9,5), (11,3)\\
16&6&4&(3,13), (5,11), (7,9), (9,7), (11,5), (13,3)\\
18&7&4&(3,15), (5,13), (7,11), (9,9), (11,7), (13,5), (15,3)\\
20&8&4&(3,17), (5,15), (7,13), (9,11), (11,9), (13,7), (15,5), (17,3)\\
22&9&5&(3,19), (5,17), (7,15), (9,13), (11,11), (13,9), (15,7), (17,5), (19,3)\\
24&10&6&(3,21), (5,19), (7,17), (9,15), (11,13), (13,11), (15,9), (17,7), (19,5), (21,3)\\
etc.&etc.&etc.&etc.\\
n&(n/2)-2&?&(3,n-3), (5,n-5),(7,n-7), ... , (n-7,7),(n-5,5),(n-3,3)\\
\hline
\end{tabular}}
\caption{Relation between even integer $n$ and number of pairs (x,y) such that x+y=n and x and y are odd integers greater than 1.}
\end{table}

To find the prime pairs (x,y) where x+y = n, we must find all the composite pairs where the $x$ or $y$ coordinate is divisible by 3. Then this value can be subtracted from |P(n)| to get the  pairs that are relatively prime to 3 (i.e. the $x$ and $y$ coordinates are not divisible by 3). The the process can be repeated to find all pairs where $x$ and $y$ are relatively prime to 5, 7, 11 etc. until we reach a prime number $p_i$, such that $p_i^2 < \sqrt{n}$ and the next higher prime $p_{i+1}^2 > \sqrt{n}$. 

\section{Number of pairs divisible by 3}

First we have to find the pairs in P(n) where either the x or y coordinate is divisible by 3.
There are two cases:

Case 1: $n$ is an even integer > 6 and $n$ is evenly divisible by 3.\\
In this case, when $x$ is evenly divisible by 3, the value of $y=n-x$ is also evenly divisible by 3. There is no instance where $x$ is divisible by 3 and $y$ is not. There are exactly $\floor{(|P(n)|-1)/3}$ pairs where both the $x$ and $y$ coordinates are evenly divisible by 3. As $n \to \infty$, this value approaches $\frac{1}{3}|P(n)|$.

For example, for $n =990$, out of the $493$ pairs, there were 165 pairs where either the $x$ or $y$ coordinate is divisible by 3. This is approximately equal to $\frac{1}{3}|P(n)| = 164\frac{1}{3}$\\

Case 2: $n$ is an even integer and $n$ is not evenly divisible by $3$. \\
In this case, there are no pairs where both the $x$ and $y$ coordinates of a pair are divisible by 3. 
There are exactly $\floor{(|P(n)| - 1)/3)}$ pairs where the $x$ coordinate is divisible by 3. Note that we have to subtract $1$ from $|P(n)|$ since we do not count the pair $(3,n-3)$. Likewise, for the $y$ coordinate, there are $\floor{(|P(n)| - 1)/3}$ pairs where the $y$ coordinate is divisible by 3 excluding $(n-3,3)$. Since there are no pairs that have both $x$ and $y$ coordinates divisible by 3, the total number of pairs where either the $x$ or $y$ coordinate is divisible by 3 is $2 \times \floor{(|P(n)|-1)/3)]}=K_3(n)$. As $n \to \infty$, this value approaches $\frac{2}{3}|P(n)|$. This value is twice the value for when $n$ is evenly divisible by 3.\\

For example, for $n =994$, out of the $495$ pairs, there are $328$ pairs such that are divisible by 3. This is approximately equal to $\frac{2}{3}|P(n)| = 330$.  Note that this is about twice as many pairs compared to $n=990$.\\

\textbf{Postulate 1: For the pairs (x,y) such that x+y=n and both x and y are odd integers, the number of pairs that are evenly divisible by 3 is relatively larger when $n$ is not evenly divisible by 3.}\\

Let $K_3(n)$ equal the number of pairs of $n$ that are divisible by 3. For a very large even $n$ not divisible by 3, the number of pairs evenly divisible by 3 approaches the following equation:
\begin{align}
K_3(n) = 2 \times \floor{(|P(n)|-1)/3)}\label{eq:eq1} \\
K_3(n) \lim_{n\to\infty} = (2/3)|P(n)|\\ 
K^*_3(n) = (2/3)|P(n)|
\end{align}
The asterisk (*) will indicate that $K^*_3(n)$ is an approximation of the actual number of pairs divisible by $n$, $K_3(n)$.

If we subtract $K^*_3(n)$ from |P(n)|, this gives us the number of pairs less than n that are not evenly divisible by 3.  As $n\to\infty$, this can be represented by the following equation.
\begin{align}
\Pi^*_3(n) &=  |P(n)| - K^*_3(n)\\
\Pi^*_3(n) &=  |P(n)| - (2/3)|P(n)| \\
\Pi^*_3(n) &=  (1/3)|P(n)|
\end{align}

\section{Number of pairs divisible by 5 and not 3}

There are two cases for calculating the number of pairs where either the $x$ or $y$ coordinate is divisible by 5:

Case 1: $n$ is an even integer and $n$ is evenly divisible by 5.\\
The reasoning is very similar a that for the case where $n$ is divisible by 3. In this case, when $x$ is evenly divisible by 5, the value of $y=n-x$ is also evenly divisible by 5. There is no instance where $x$ is divisible by 5 and $y$ is not. There are exactly $\floor{(P(n) - 2)/5}$ pairs where both the $x$ and $y$ coordinates are divisible by 5. As $n \to \infty$, this value approaches $\frac{1}{5}|P(n)|$.\\

For example, for $n =990$, out of the 493 pairs, there are $99$ pairs where both $x$ and $y$ are divisible by 5. This is approximately equal to  $\frac{1}{5}P(n) = 98.6$. Of these $99$ pairs, $33$ pairs of these pairs are also divisible by 3.\\

Case 2: $n$ is even and $n$ is not evenly divisible by 5.\\
In this case, there are no pairs where both the $x$ and $y$ coordinates of a pair are divisible by 5. 
There are exactly $\floor{(|P(n)|- 2)/5}$ pairs where the $x$ coordinate is divisible by 5. Note that we have to subtract 2 from $|P(n)|$ since we do not count the first two pairs $(3,n-3)$ and $(5,n-5)$.
Likewise for the $y$ coordinate, there are $\floor{(|P(n)| - 2)/5)}$ pairs where the $y$ coordinate is divisible by 5.
Since there are no pairs that have both $x$ and $y$ coordinates divisible by 5, the total number of pairs where either the $x$ or the $y$ coordinate is divisible by 5 is 
\begin{align}
k_5(n) = 2 \times \floor{(|P(n)| - 2)/5)} \label{eq:eq2}
\end{align}
As $n \to \infty$, this value approaches 
\begin{align}
k_5(n) \lim\to\infty = \frac{2}{5}|P(n)|
\end{align}

For example, for $n=994$ there are $196$ pairs that are divisible by 5 compared to when $n=990$ there are only 99 pairs divisible by 5. There are about twice as many pairs divisible by 5 when $n$ is not divisible by 5.

\textbf{Postulate 2: For the pairs (x,y) such that x+y=n and both x and y are odd integers, the number of pairs that are evenly divisible by 5 is relatively larger when $n$ is not evenly divisible by 5.}\\

Of the pairs where the $x$ coordintate is divisible by 5, every third pair is divisible by 3. Likewise for the $y$ coordinate. So for large values of $n$, about 2/3rds of the pairs are divisible by 3 and must be subracted from the pairs divisible by 5. Subtracting 2/3rds of the pairs is like multiplying by 1/3. So for very large $n$ not divisible by 5 or 3, the number of pairs divisible by 5 and not 3 approaches the following equation:\\
\begin{align}
K_5(n) \lim_{n\to\infty} = (1/3)(2/5)|P(n)|\\
K^*_5(n)  = (1/3)(2/5)|P(n)|
\end{align}

Subtracting $K^*_5(n)$ from $\Pi^*_3(n)$ gives us the number of pairs not divisible by 3 or 5. As $n \to \infty$, this can be represented by the following equation.
\begin{align}
\Pi^*_5(n) &=  \Pi^*_3(n) - K^*_5(n)\\
\Pi^*_5(n) &=  (1/3)|P(n)| -(1/3) (2/5)|P(n)|\\
\Pi^*_5(n) &= (1/3)(3/5)|P(n)| \label{eq:pi5}
\end{align}

For example, for $n=994$, of the $196$ pairs are divisible by 5, there are 132 pairs (about two-thirds) that are divisible by 3. This is because for every $x$ coordinate divisible by 5 starting with (15,979), every third pair is also divisible by 3 (yellow numbers) and starting with (25,969), every third pair the $y$ coordinate is divisible by 3 (orange numbers). There are no cases where both $x$ and $y$ are divisible by 3.\\
(\colorbox{yellow}{15},979),(25,\colorbox{orange}{969}),(35,959),(\colorbox{yellow}{45},949),(55,\colorbox{orange}{939}),(65,929),(\colorbox{yellow}{75},919),(85,\colorbox{orange}{909}),\\(95,899),(\colorbox{yellow}{105},889),(115,\colorbox{orange}{879})...\\

So about one third (64 pairs) of the 196 pairs divisible by 5 are not divisible by 3.\\

\section{Number of pairs divisible by 7}

There are two cases for calculating the number of pairs where either the $x$ or $y$ coordinate is divisible by 7:

Case 1: $n$ is an even integer and $n$ is evenly divisible by 7.\\
In this case, when $x$ is evenly divisible by 7, the value of $y=n-x$ is also evenly divisible by 7. There is no instance where $x$ is divisible by 7 and $y$ is not. There are exactly $\floor{(|P(n)| - 3)/7}$ pairs where both the $x$ and $y$ coordinates are divisible by 7. As $n \to \infty$, this value approaches $\frac{1}{7}|P(n)|$.\\

Case 2: $n$ is an even integer and $n$ is not evenly divisible by 7.\\
In this case, there are no pairs where both the $x$ and $y$ coordinates of a pair are divisible by 7. 
There are exactly $\floor{(|P(n)|- 3)/7}$ pairs where the $x$ coordinate is divisible by 7. Note that we have to subtract 3 from $|P(n)|$ since we do not count the first three pairs $(3,n-3)$, $(5,n-5)$ and  $(7,n-7)$.
Likewise for the $y$ coordinate, there are $\floor{(|P(n)| - 3)/7)}$ pairs where the $y$ coordinate is divisible by 7.
So the total number of pairs where either the $x$ or the $y$ coordinate is divisible by 7 is $2 \times \floor{(|P(n)| - 3)/7)}$. As $n \to \infty$, this value approaches $\frac{2}{5}|P(n)|$.\\

\textbf{Postulate 3: For the pairs (x,y) such that x+y=n and both x and y are odd integers, the number of pairs that are evenly divisible by 7 is relatively larger when $n$ is not evenly divisible by 7.}\\

About 2/5ths of the pairs divisible by 7 are also divisible by 5. Subtracting 2/5ths of the pairs is like multiplying by 3/5. 

About 2/3rds of the remaining pairs are also divisible by 3. Subtracting 2/3rds of the pairs is like multiplying by 1/3. 

For very large $n$, the number of pairs divisible by 7 and not 5 or 3 approaches the following equation:\\
\begin{align}
K_7(n) \lim_{n\to\infty} = (1/3)(3/5)(2/7)|P(n)|\\
K^*_7(n) = (1/3)(3/5)(2/7)|P(n)|
\end{align}

Subtracting $K_7(n)$ from $Pi^*_5(n)$ gives us the number of pairs not divisible by 3, 5 or 7. As $n \to \infty$, this can be represented by the following equation:

\begin{align}
\Pi^*_7(n) &=  \Pi^*_5(n) - K^*_5(n)\\
\Pi^*_7(n) &=  (1/3)(3/5)|P(n)| -(1/3)(3/5)(2/7)|P(n)|\\
\Pi^*_7(n) &= (1/3)(3/5)(5/7)|P(n)|
\end{align}

\section{Number of pairs divisible by a prime > 7}

Using the same technique, we can calculate the number of pairs that are divisible by any prime number $p_i$ and for each prime number, there are two cases.
Case 1: $n$ is an even integer and $n$ is evenly divisible by $p_i$.\\
In this case, when $x$ is evenly divisible by $p$, the value of $y=n-x$ is also evenly divisible by $p$. There is no instance where $x$ is divisible by $p$ and $y$ is not. There are exactly $\floor{(|P(n)| - (p-1)/2)/p}$ pairs where both the $x$ and $y$ coordinates are divisible by $p$. As $n \to \infty$, this value approaches $\frac{1}{p}|P(n)|$.\\

Case 2: $n$ is an even integer and $n$ is not evenly divisible by $p_i$.\\
In this case, there are exactly $2\times\floor{((|P(n)|- (p-1)/2)/p}$ pairs where $x$ or $y$ is divisible by $p$. We have to subtract $(p-1)/2$ from $|P(n)|$ since we do not count the pairs $(3,n-3), (5,n-5), (7,n-7)$,...,$(p,n-p)$. We have to multiply by 2 since the $x$ coordinate and the $y$ coordinate cannot both be divisible by $p$. 

\textbf{Postulate 4: For pairs $(x,y)$ such that $x+y=n$ and both $x$ and $y$ are odd integers, the number of pairs that are evenly divisible by $p$ is relatively larger when $n$ is not evenly divisible by $p$.}\\

Combining posutates 1-4, the maximum number of pairs $(x,y)$ such that $x+y=n$ and both $x$ and $y$ are odd integers, where $x$ or $y$ is a composite number is relatively largest when $n$ is not evenly divisible by any prime number less than $l(\sqrt{n})$. In other words, when $n=2^i$ or $n=2^ip$ where $i$ is an integer greater than 0 and $p$ is a prime number greater than $l(\sqrt{n})$, the number of prime pairs of $n$ is fewest. This is why the curve where $n$ is not divisible by a prime less than $\sqrt{n}$, the orange line in Figure \ref{fig:Graphs}B, is the lower bound and the curves where $n$ has one or more prime factors greater than 2 have relatively more prime pairs.\\

Using this technique, the values for $\Pi^*_p(n)$ can be determined to be as follows:
\begin{equation}
\Pi^*_p(n) = |P(n)| \times (1/3)(3/5)(5/7)(9/11) ... \frac{l(\sqrt{n})-2)}{l(\sqrt{n})}
\end{equation}

Note that we only need to check if the $x$ or $y$ coordinate is divisible by prime numbers less than or equal to  $p_{i} = l(\sqrt{n})$. For any prime $p_{i+1} > l(\sqrt{n})$, the first pair where the $x$ coordinate wouuld be divisible by $p_{i+1}$ and not already counted is $p_{i+1}^2$ which is greater than $n$.

The total number of prime pairs of $n$ where $n$ is not divisible by a prime number less than $l(\sqrt{n})$,  can be approximated by the following equation:
\begin{equation}
\Pi^*(n) = |P(n)|\prod_{(p=3)}^{(p=l(\sqrt{n}))} \frac{p-2}{p}
\end{equation}
where the product is over prime numbers.

Let us define the function $W(n)$, which represents the fraction of pairs of $P(n)$ where the $x$ and $y$ coordinates are prime.\\
\begin{equation}
W(n) = \prod_{(p=3)}^{(p=l(\sqrt{n}))} \frac{p-2}{p}
\end{equation}
where the product is over prime numbers.

This simplifies the equation for the number of prime pairs (x,y) where x+y=n to the following:
\begin{equation}
\Pi^*(n) = |P(n)|W(n)\label{eq:pistar}
\end{equation}

\begin{figure}
  \includegraphics[width=\linewidth]{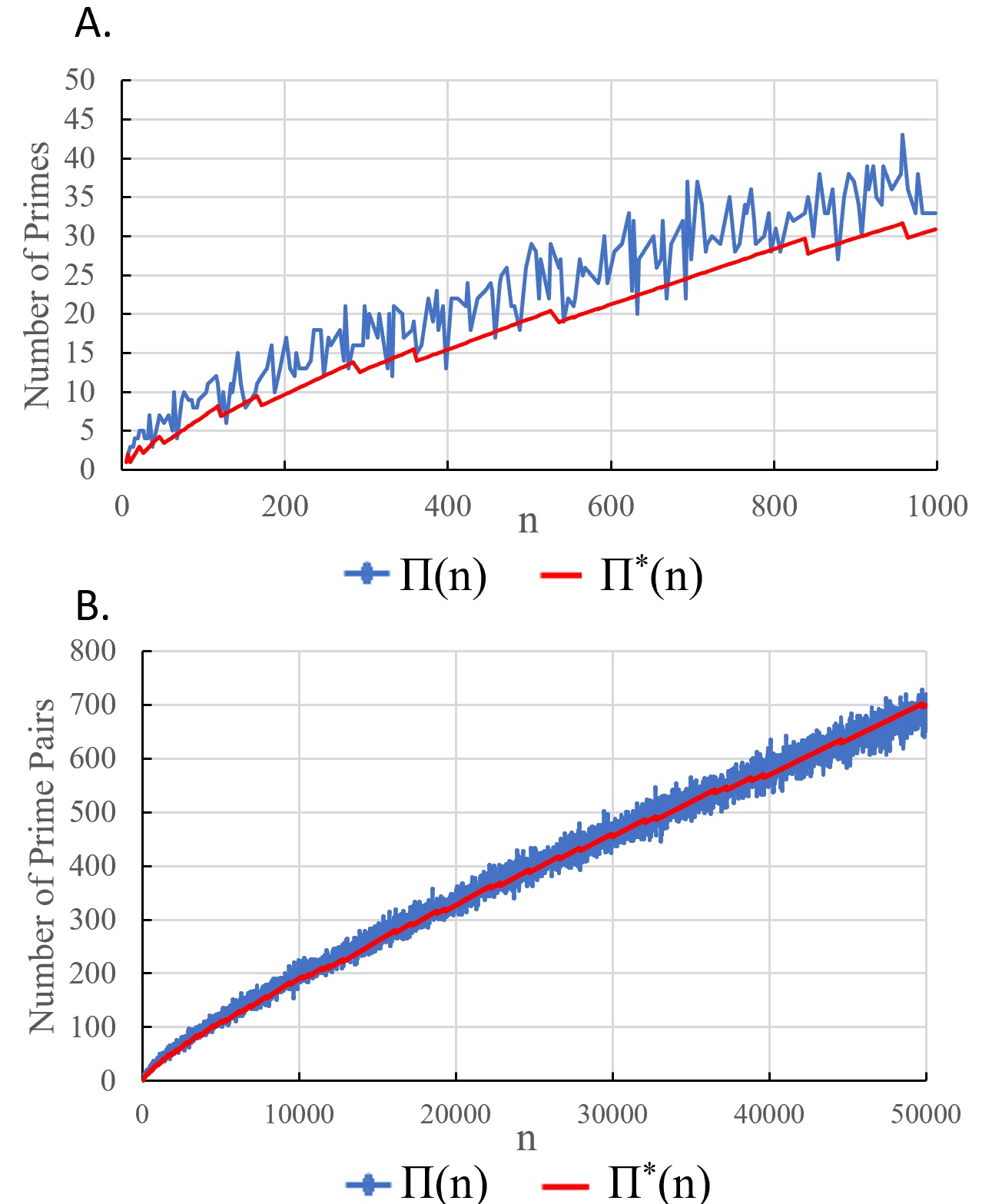}
  \caption{$\Pi^*(n)=P(n)W(n)$ accurately predicts the number of prime pairs for n not divisible by a prime number less than $\sqrt{n}$. $\Pi^*(n)$ modestly underestimates the number of prime pairs (A.) for values of $n$ less than 1000, but as $n$ increases (B.) the number of prime pairs approaches $P(n)W(n)$.}
  \label{fig:Graphs2}
\end{figure}

$\Pi^*(n)$ approximates the number of prime pairsof for $n$ not evenly divisible by a prime number less than $\sqrt{n}$. To confirm that no error was made in the derivation of this equation and to visualize the accuracy of the equations, the values for $\Pi^*(n)$ were plotted against values of $n$ that are not evenly divisible by a prime number less than $\sqrt{n}$ between 6 and 50,000 and compared them to the actual number of prime pairs $\Pi(n)$. As can be seen in Figure \ref{fig:Graphs2}A, the curve for $\Pi^*(n)$ (red curve) modestly underestimates the actual number of prime pairs (blue curve) for $n$ less than 1,000. But as $n$ gets large, (Figure \ref{fig:Graphs2}B) $\Pi^*(n)$ does indeed approach the actual number of prime pairs of $n$. Also note the jagged saw-tooth appearance of the curve for $\Pi^*(n)$. This occurs because the $W(n)$ function increases stepwise every time the value of $n$ changes from $p^2-1$ to $p^2+1$ as depicted in Figure \ref{fig:zigzag}.

\section{Proof of Goldbach's Conjecture}
In order to perform an inductive proof, we need to understand how the $W(n)$ function changes as $n$ increases.
Let us look at the values of $W(n)$ where $n=p_i^2+1$ since $W(n)$ does not change for values of $n = p_i^2+1$ to $n=p_{i+1}^2-1$. This is done to determine the relationship between $W(p_i^2+1)$ and $W(p_{i+1}^2+1)$ so that a proof by induction can be performed.\\

Notice that $W(p_{i+1}^2+1)$ can be expressed as:\\
\begin{equation}
W(p_{i+1}^2+1) = \left(\frac{p_{i+1}-2}{p_{i+1}}\right)W(p_i^2+1)
\end{equation}

Also notice that $W(p_{i+1}^2-1)$ can be expressed as:\\
\begin{equation}
W(p_{i+1}^2-1) = \left(\frac{p_i-2}{p_i}\right)W(p_i^2-1)
\end{equation}

To prove Goldbach's conjecture, we must prove that $\Pi^*(n) \geq 1$ for all even $n$ greater than 4.
There are three cases that we will look at and for all cases, we will prove that $\Pi^*(n) \geq 1$.

\begin{figure}
  \includegraphics[scale=0.9]{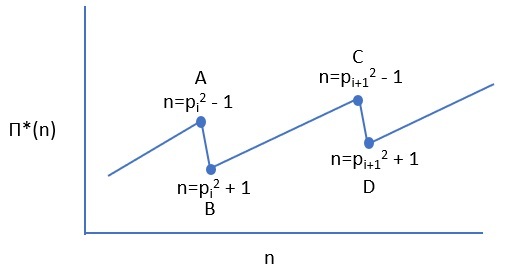}
  \caption{This figure depicts the jagged saw-tooth appearance of $\Pi^*(n)$ versus $n$. The drops occur when $n$ goes from $p_i^2-1$ to $p_i^2+1$.}
  \label{fig:zigzag}
\end{figure}

Case 1 is where $p_i^2 + 1 \leq n < p_{i+1}^2 -1$. We will prove that for all $n$ such that  $p_i^2 + 1 \leq n < p_{i+1}^2 -1$, $\Pi(n) \geq 1$ and $\Pi(n+2) \ge \Pi(n)$. This would be the region between points B and C in Figure \ref{fig:zigzag}. 

Case 2 is where $n=p^2-1$. For $n_i = p_i^2 -1$ and $n_{i+1} = p_{i+1}^2 -1$, we will prove that for all $n_i$, 
$\Pi(n_i) \ge 1$ and that $\Pi(n_{i+1}) \geq \Pi(n_i)$ for all $n_i$. In other words, we will show that point C is always greater than or equal to point A in Figure \ref{fig:zigzag}.

Case 3 is where $n=p^2+1$. For $n_i = p_i^2 +1$ and $n_{i+1} = p_{i+1}^2 +1$, we will prove that for all $n_i$,
$\Pi(n_i) \ge 1$ and that $\Pi(n_{i+1}) \geq \Pi(n_i)$ for all $n_i$. In other words, we will show that point D is always greater than or equal to point B in Figure \ref{fig:zigzag}.\\

Case 1: $p_i^2 + 1 \leq n < p_{i+1}^2 -1$. The region between points B and C in Figure \ref{fig:zigzag}\\

For values of $n$, such that $p_i^2 + 1 \leq n < p_{i+1}^2 -1$, the value of $W(n)$ does not change.
\begin{align}
\Pi^*(n) &= |P(n)|W(n)\\
\Pi^*(n) &= (\frac{n}{2} -2)W(n)
\end{align}
Since $W(n)$ is constant, this equation increases linearly with increasing $n$.\\

Case 2: For $n = p^2 - 1$, $\Pi(n) > 1$ and $\Pi(p_{i+1}^2-1) > \Pi(p_i^2-1)$. Point C is greater than point A in Figure \ref{fig:zigzag}. 
To prove this case, we will also use proof by induction.\\

Base case: For the base case $p = 3$ and $n=3^2 - 1= 8$.\\
\begin{align}
\Pi^*(8) &= |P(8)|W(8)\\
\Pi^*(8) &= ((8/2)-2)(1)\\
\Pi^*(8) &= 2 > 1.
\end{align}

Inductive step: We will assume that $\Pi^*(p_i^2-1) \geq 1$, and prove that $\Pi^*(p_{i+1}^2-1) \geq 1$.\\
Assumption: $\Pi^*(p_i^2-1) \geq 1$.\\
Prove: $\Pi^*(p_{i+1}^2-1) \geq 1$.\\
\begin{align}
\Pi^*(p_{i+1}^2-1)&=P(p_{i+1}^2-1)W(p_{i+1}^2-1)\\
\Pi^*(p_{i+1}^2-1) &= P(p_{i+1}^2-1)\left(\frac{p_i-2}{p_i}\right)W(p_i^2-1)\\
\Pi^*(p_{i+1}^2-1) &= P(p_{i+1}^2-1)\left(\frac{p_i-2}{p_i}\right)\left(\frac{P(p_i^2-1)}{P(p_i^2-1)}\right)W(p_i^2-1)\\
\Pi^*(p_{i+1}^2-1) &= \left(\frac{P(p_{i+1}^2-1)}{P(p_i^2-1)}\right) \left(\frac{p_i-2}{p_i}\right)\Pi^*(p_i^2-1)\\
\Pi^*(p_{i+1}^2-1) &= \left(\frac{(p_{i+1}^2-1)/2-2}{(p_i^2-1)/2-2}\right) \left(\frac{p_i-2}{p_i}\right)\Pi^*(p_i^2-1)\\
\Pi^*(p_{i+1}^2-1) &= \left(\frac{(p_{i+1}^2-5}{p_i^2-5}\right) \left(\frac{p_i-2}{p_i}\right)\Pi^*(p_i^2-1)
\end{align}
Since we assumed that $\Pi^*(p_i^2-1) \geq 1$, we need to show that the term
$\left(\frac{p_{i+1}^2-5}{p_i^2-5}\right)\left(\frac{p_i-2}{p_i}\right) \geq1$. 

\begin{align}
\left(\frac{p_{i+1}^2-5}{p_i^2-5}\right)\left(\frac{p_i-2}{p_i}\right) &\geq1\\
(p_{i+1}^2-5)(p_i-2) &\geq (p_i^2-5)(p_i)\\
p_{i+1}^2p_i - 2p_{i+1}^2 - 5p_i +10 &\geq (p_i^3-5p_i)\\
p_{i+1}^2p_i - 2p_{i+1}^2 +10 &\geq p_i^3
\end{align}
To minimize the LHS, the lowest value $p_{i+1}$ can be is $p_i+2$. So substituting gives:
\begin{align}
(p_i+2)^2p_i - 2(p_i+2)^2 +10 &\geq p_i^3\\
(p_i^2 + 4p_i+ 4)p_i - 2(p_i^2 +4p_i+4) +10 &\geq p_i^3\\
p_i^3 + 4p_i^2+ 4p_i - 2p_i^2 -8p_i-8 +10 &\geq p_i^3\\
4p_i^2+ 4p_i - 2p_i^2 -8p_i+2 &\geq 0\\
2p_i^2 - 4p_i+2 &\geq 0\\
2(p_i-1)^2 &\geq 0
\end{align}
Since the inequality $p_i-1 \geq 0$ holds true for all $p_i$, this proves that the number of prime pairs less than $n$ increases as $n$ goes from  $n=p_i^2-1$ to $n=p_{i+1}^2-1$.\\

Case 3: For $n = p^2 + 1$, $\Pi(n) > 1$ and $\Pi(p_{i+1}^2+1) > \Pi(p_i^2+1)$. Point D is greater than point B in Figure \ref{fig:zigzag}. 
To prove this case, we will also use proof by induction. The base case $p = 3$ and $n=p^2 + 1= 10$.\\
\begin{align}
\Pi^*(10) = P(10)W(10)\\
\Pi^*(10) = ((10/2)-2)(1/3)\\
\Pi^*(10) = (3)(1/3) = 1.
\end{align}
Now assuming that $\Pi^*(p_i^2+1) \geq 1$, we must prove that $\Pi^*(p_{i+1}^2+1) \geq 1$.\\
Assumption: $\Pi^*(p_i^2+1) \geq 1$.\\
Prove: $\Pi^*(p_{i+1}^2+1) \geq 1$.
\begin{align}
\Pi^*(p_{i+1}^2+1)&=P(p_{i+1}^2+1)W(p_{i+1}^2+1)\\
\Pi^*(p_{i+1}^2+1) &= P(p_{i+1}^2+1)\left(\frac{p_{i+1}-2}{p_{i+1}}\right)W(p_i^2+1)\\
\Pi^*(p_{i+1}^2+1) &= P(p_{i+1}^2+1)\left(\frac{p_{i+1}-2}{p_{i+1}}\right)\left(\frac{P(p_i^2+1)}{P(p_i^2+1)}\right)W(p_i^2+1)\\
\Pi^*(p_{i+1}^2+1) &= \left(\frac{P(p_{i+1}^2+1)}{P(p_i^2+1)}\right) \left(\frac{p_i-2}{p_i}\right)\Pi^*(p_i^2+1)\\
\Pi^*(p_{i+1}^2+1) &= \left(\frac{(p_{i+1}^2+1)/2-2}{(p_i^2+1)/2-2}\right) \left(\frac{p_i-2}{p_i}\right)\Pi^*(p_i^2+1)\\
\Pi^*(p_{i+1}^2+1) &= \left(\frac{(p_{i+1}^2-3}{p_i^2-3}\right) \left(\frac{p_i-2}{p_i}\right)\Pi^*(p_i^2+1)
\end{align}
Since we assumed that $\Pi^*(p_i+1) \geq 1$, we need to show that the item 
$\left(\frac{p_{i+1}^2-3}{p_i^2-3}\right)\left(\frac{p_i-2}{p_i}\right) \geq1$.

\begin{align}
\left(\frac{p_{i+1}^2-3}{p_i^2-3}\right)\left(\frac{p_i-2}{p_i}\right) \geq1\\
(p_{i+1}^2-3)(p_i-2) \geq (p_i^2-3)p_i\\
p_{i+1}^2p_i-2p_{i+1}^2 -3p_i+6 \geq p_i^3-3p_i\\
p_{i+1}^2p_i-2p_{i+1}^2+6 \geq p_i^3
\end{align}
To minimize the LHS, the lowest value $p_{i+1}$ can be is $p_i+2$. So substituting gives:
\begin{align}
(p_i+2)^2p_i - 2(p_i+2)^2 +6 &\geq p_i^3\\
(p_i^2 + 4p_i+ 4)p_i - 2(p_i^2 +4p_i+4) +6 &\geq p_i^3\\
p_i^3 + 4p_i^2+ 4p_i - 2p_i^2 -8p_i-8 +6 &\geq p_i^3\\
4p_i^2+ 4p_i - 2p_i^2 -8p_i-2 &\geq 0\\
2p_i^2 - 4p_i-2 &\geq 0\\
2(p_i-1)^2 -4&\geq 0
\end{align}
Since the inequality $2(p_i-1)^2 -4 \geq 0$ holds true, this proves that the number of prime pairs less than $n$ increases as $n$ goes from  $n=p_i^2+1$ to $n=p_{i+1}^2+1$.\\

By combining these three cases, this proves that for $n=2^i$ or $n=2^ip$ where $i$ is an integer greater than 0 and $p$ is a prime number greater than $\sqrt{n}$, there is at least 1 prime pair $(p_x,p_y)$ such that $p_x+p_y = n$ and the number of prime pairs increases with increasing $n$. Since $n=2^i$ or $n=2^ip$ is the lower bound on the number of prime pairs of $n$, this proves that the Goldbach Conjecture holds true for all even integers for sufficiently large $n$. This also proves that the number of prime pairs increases as $n$ increases.

\section{Error Analysis}

We have shown that $\Pi^*(n) = P(n)W(n)$ is an approximation of the number of prime pairs of $n$ where $n$ is not evenly divisible by a prime number less than $\sqrt{n}$. As $n$ increases, this approximation approaches the true number of prime pairs of $n$, $\Pi(n)$. Using this approximation, we have proven that the number of prime pairs of $n$ increases as $n$ increases, thus proving the Golbach's conjecture is true for sufficiently large $n$. The next question is, how large is a sufficiently large value of $n$?\\

To do this, we have to find the error in the $\Pi^*(n)$ function. Let us examine how the $\Pi^*(n)$ was derived to determine how much $\Pi^*(n)$ deviates from $\Pi(n)$\\

We will first calculate the maximum error in $\Pi^*_3(n)$. We will use a lowercase $k_p(n)$ to differentiate from $K_p(n)$. The lowercase $k_p(n)$ indicates the number of pairs of $n$ that are divisible by $p$. We will also use $\bar k_p(n)$ to indicate the number of pairs that are not divisible by $p$, i.e. $\bar k_p(n) = |P(n)| - k_p(n)$.

The exact number of pairs of $n$ that are divisible by 3 was defined by equation \ref{eq:eq1}:
\begin{align}
k_3(n) = 2\floor{(|P(n)|-1)/3)}
\end{align}
Therefore, the exact number of pairs of $n$ that are not divisible by 3, $\bar k_3(n)$, is defined by the equation:
\begin{align}
\bar k_3(n) &=  |P(n)| - k_3(n) \\
\end{align}

As $n$ gets large, $k_3(n)$ was approximated by the following equation:
\begin{align}
k^*_3(n) = (2/3)|P(n)|
\end{align}
Therefore, the approximation of the number of pairs of $n$ that are not divisible by 3, $\bar k^*_3(n)$, is defined by the equation:
\begin{align}
\bar k^*_3(n) &=  |P(n)| - k^*_3(n)
\end{align}

The difference between $\bar k_3(n)$ and $\bar k^*_3(n)$ is the error, $e_3(n)$ and can be defined as follows:
\begin{align}
e_3(n) &= \bar k_3(n) - \bar k^*_3(n) \\
e_3(n) &= |P(n)| - k_3(n) - (|P(n)| - k^*_3(n)) \\
e_3(n) &= k^*_3(n) - k_3(n)) \\
e_3(n) &= (2/3)|P(n)| - 2\floor{(|P(n)|-1)/3)} \label{eq:e3}
\end{align}

\begin{figure}
  \includegraphics[width=\linewidth]{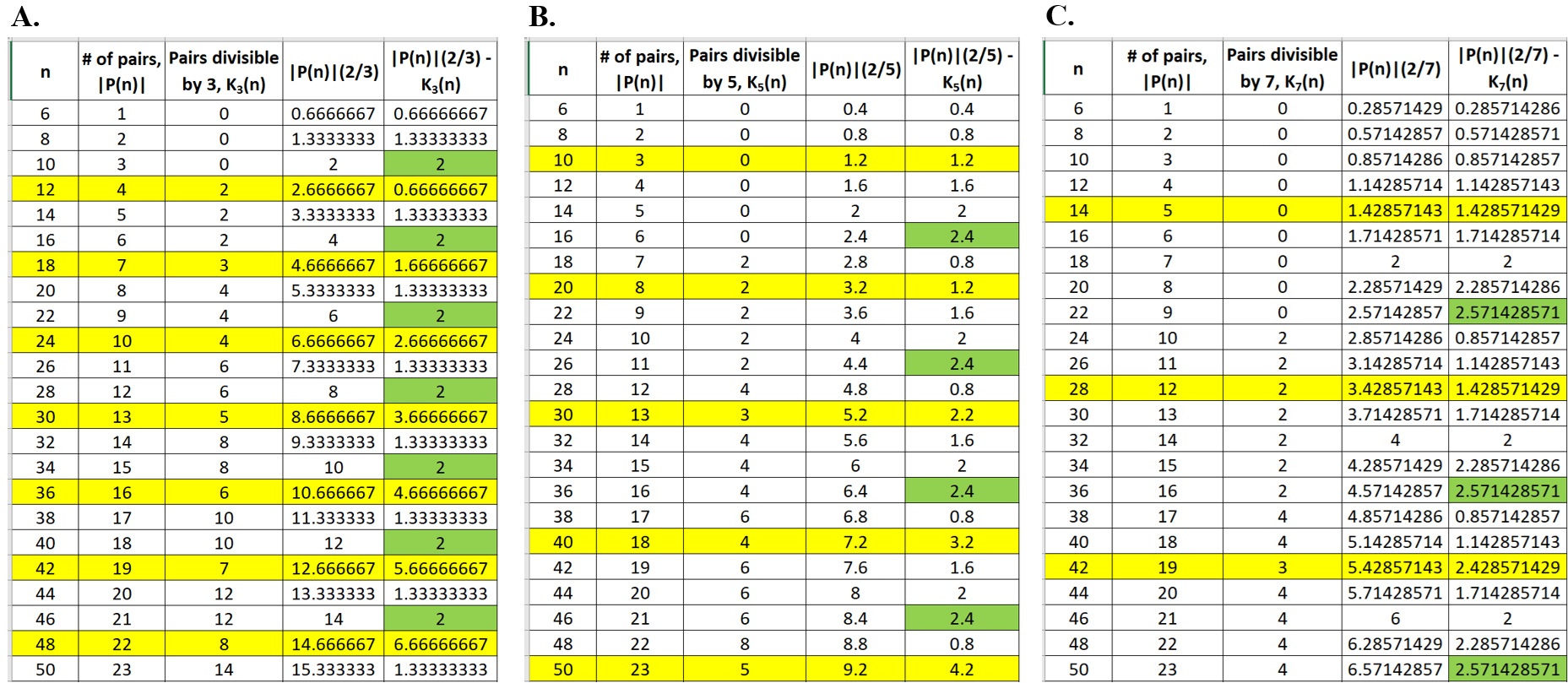}
  \caption{The difference between the actual number of pairs divisible by 3, 5 or 7 and |P(n)|(2/3), |P(n)|(3/5) and |P(n)|(5/7) respectively. The maximum errors are highlighted in green. The values of n that are evenly divisible by 3, 5 or 7 are highlighted in yellow and are excluded.}
  \label{fig:Pairs}
\end{figure}

The maximum value of $e_3(n)$ for even $n$ occurs when $|P(n)| = 3+3i$ where $i \geq 0$. This is because $|P(n)| =3$ is the largest value of $|P(n)|$ where the floor function of $e_3$ rounds down to 0. Plugging numbers into equation \ref{eq:e3} confirms that the maximum value of $e_3(n)$ is 2 as shown in Figure \ref{fig:Pairs}A.

\begin{align}
e_3 = 2
\end{align}
Note that $e_3$ is a constant and will be greater than or equal to 2 since we are excluding the pairs (3,n-3) and (n-3,3).

Next, we can calculate the maximum error in $k^*_5(n)$. The exact number of pairs of $n$ that are divisible by 5 was defined by equation \ref{eq:eq2}:
\begin{align}
k_5(n) = 2\floor{(|P(n)|-2)/5)}
\end{align}
Therefore, the exact number of pairs of $n$ that are not divisible by 5, $\bar k_5(n)$, is defined by the equation:
\begin{align}
\bar k_5(n) &=  |P(n)| - k_5(n) \\
\end{align}

As $n$ gets large, $k_3(n)$ was approximated by the following equation:
\begin{align}
k^*_5(n) = (2/5)|P(n)|
\end{align}
Therefore, the approximation of the number of pairs of $n$ that are not divisible by 5, $\bar k^*_3(n)$, is defined by the equation:
\begin{align}
\bar k^*_5(n) &=  |P(n)| - k^*_5(n)
\end{align}

The difference between $\bar k_5(n)$ and $\bar k^*_3(n)$ is the error, $e_5(n)$ and can be defined as follows:
\begin{align}
e_5(n) &= \bar k_5(n) - \bar k^*_5(n) \\
e_5(n) &= |P(n)| - k_5(n) - (|P(n)| - k^*_5(n)) \\
e_5(n) &= k^*_5(n) - k_5(n)) \\
e_5(n) &= (2/5)|P(n)| - 2\floor{(|P(n)|-2)/5)} \label{eq:e5}
\end{align}
The maximum value of $e_5(n)$ for even $n$ occurs when $|P(n)| = 6+5i$ where $i \geq 0$. This is because $|P(n)| =6$ is the largest value of $|P(n)|$ where the floor function of $e_5$ rounds down to 0. Plugging numbers into equation \ref{eq:e5} confirms that the maximum value of $e_5(n)$ is 2.4 as shown in Figure \ref{fig:Pairs}B.

\begin{align}
e_5 = 2.4
\end{align}
Note that $e_5$ is a constant.

Next, we need to determine the maximum error in $\Pi^*_5(n)$. To do this, let us define the following functions:
\begin{center}
$a(n)$ = number of pairs divisible by 3 excluding (3,n-3) and (n-3,3). \\
$b(n)$ = number of pairs divisible by 5 excluding (5,n-5) and (n-5,5).\\
$c(n)$ = number of pairs divisible by both 3 and 5\\
\end{center}
In this case $a(n) = k_3(n)$ and $b(n) = k_5(n)$. The function $c(n)$ would include all pairs divisible by 15 plus pairs where the x coordinate is divisible by 3 and the y coordinate is divisible by 5 plus pairs where the x coordinate is divisible by 5 and the y coordinate is divisible by 3. The function c(n) can be defined as follows:
\begin{align}
c(n) = 2\floor{(P(n)+8)/15} + \floor{(P(n)+u)/15} + \floor{(P(n)+v)/15}
\end{align}
where $u$ and $v$ are integers that depend on the first occurence of a pair of $n$ where $x$ is divisible by 3 and $y$ is divisible by 5 and where $x$ is divisible by 5 and $y$ is divisible by 3 respectively.

The inclusion/exclusion principle states that given two set A and B, the number of elements in A union B is $|A \cup B| = |A| + |B| - |A \cap B|$. Also note that $|A \cup B| \leq |A| + |B|$. Using the inclusion/exclusion principle, we can define number of pairs of $n$ that are evenly divisible by 3 or 5, $K_5(n)$, as follows:
\begin{align}
K_5(n) = a(n) + b(n) - c(n)
\end{align}
Subtracting $K_5(n)$ from $|P(n)|$ gives us $\Pi_5(n)$:
\begin{align}
\Pi_5(n) &= |P(n)| - K_5(n)\\
\Pi_5(n) &= |P(n)| - a(n) - b(n) + c(n)\\
\Pi_5(n) &= |P(n)| - 2\floor{(P(n)-1)/3}  - 2\floor{(P(n)-2)/5} + c(n)
\end{align}

We also know from equation \ref{eq:pi5}, that 
\begin{align}
\Pi_5^*(n) &= (1/3)(1/5)|P(n)|\\
\Pi_5^*(n) &= P(n) - (2/3)P(n) - (2/5)P(n) + (4/15)P(n)
\end{align}

Now we can calculate the error in $\Pi_5^*(n)$. Let us call that error $e_{max5}(n)$.
\begin{align}
e_{max5}(n) &= \Pi_5(n) - \Pi_5^*(n)\\
e_{max5}(n) &= P(n) - 2\floor{(P(n)-1)/3} - 2\floor{(P(n)-2)/5} + c(n)\notag\\
&- (P(n) - (2/3)P(n) -  (2/5)P(n) + (4/15)P(n))\\
e_{max5}(n) &= (2/3)P(n) - 2\floor{(P(n)-1)/3} + (2/5)P(n) - 2\floor{(P(n)-2)/5}\notag\\
&- ((4/15)P(n) - c(n))
\end{align}
Since $e_3(n)$ and $e_5(n)$ were derived as follows:
\begin{align}
e_3(n) &= (2/3)|P(n)| - 2\floor{(|P(n)|-1)/3}\\
e_5(n) &= (2/5)|P(n)| - 2\floor{(|P(n)|-2)/5}
\end{align}
We can substitute for $e_3(n)$ and $e_5(n)$ to get:
\begin{align}
e_{max5}(n) = e3 + e5  - ((4/15)P(n) - c(n))
\end{align}
As we add errors $e_7$, $e_{11}$, $e_{13}$, etc. the inclusion/exclusion principle starts getting exponentially complex since we have to subtract overlapping values and add back in other overlapping values. But since $|A \cup B| \leq |A| + |B|$, we can simplify the equation as follows:
\begin{align}
e_{max5}(n) \leq e3 + e5 
\end{align}

The general formula for $e_p$ is as follows:
\begin{align}
e_p(n) &= \left(\frac{2}{p}\right)|P(n)| - 2\left\lfloor\frac{(|P(n)|-\frac{(p-1)}{2})}{p}\right\rfloor\label{eq:epn}
\end{align}

The maximum value of $e_p(n)$ occurs at the value of $|P(n)|$ where the floor function rounds down to 0. 
This occurs when 
\begin{align}
|P(n)|-\frac{(p-1)}{2} &= p-1\\
|P(n)| &= p +\frac{(p-1)}{2}-1\\
|P(n)| &= \frac{(2p +(p-1)-2)}{2}\\
|P(n)| &= (3/2)(p -1)
\end{align}

By plugging $|P(n)| = (3/2)(p -1)$ into equation \ref{eq:epn} we get the maximum error for $e_p$ as
\begin{align}
e_p &= (2/p)(3/2)(p -1) - 0\\
e_p &= (3p-3)/p\label{eq:ep}
\end{align}

Using equation \ref{eq:ep}, we can find the maximum errors for every value of $e_p$.

\begin{align}
e_3 &= 6/3 =2\\
e_5 &= 12/5 = 2.4\\
e_7 &= 18/7 = 2.571429\\
e_{11} &= 30/11 = 2.727273\\
e_{13} &= 36/13 = 2.615385
\end{align}

If we sum up all these errors up to $e_{l(\sqrt{n})}$, this will give us the maximum error between $\Pi^*(n)$ and the actual prime pairs of n, $\Pi(n)$, for $n$ not divisible by a prime number less than $\sqrt{n}$.

\begin{align}
e_{max}(n) = \sum_{p=3}^{l(\sqrt{n})}e_p\\
e_{max}(n) = \sum_{p=3}^{l(\sqrt{n})}(3p-3)/p\\
e_{max}(n) \leq 3\pi(\sqrt{n})
\end{align}

\begin{figure}
  \includegraphics[width=\linewidth]{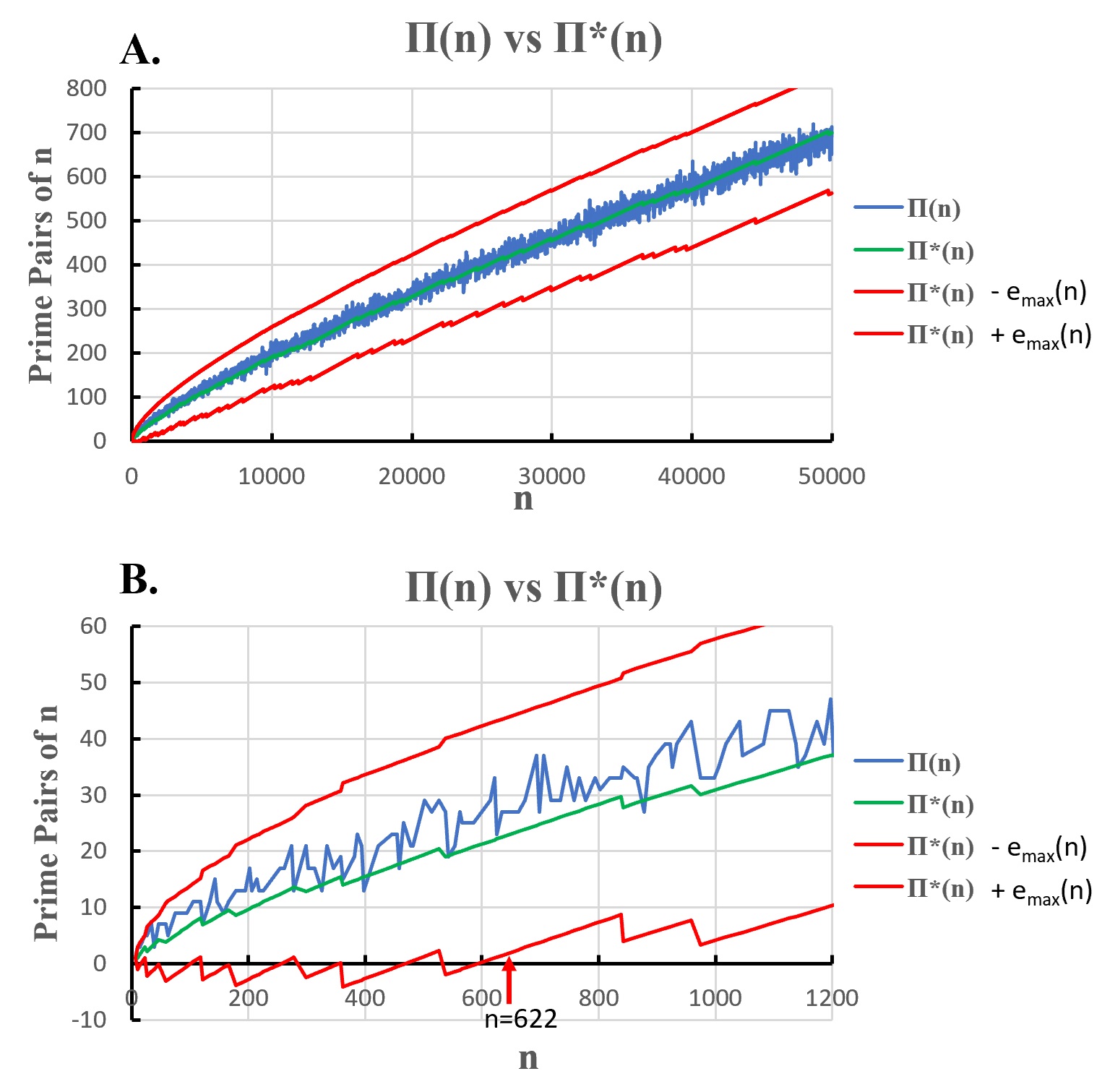}
  \caption{Graph of $\Pi(n)$ and $\Pi^*(n)$ vs. $n$ for $n$ not divisible by a prime number less than $\sqrt{n}$ shows that (A.) $\Pi(n)$ fits comfortably within $\Pi^*(n) \pm e_{max}(n)$ for all values of $n$ up to 50000. (B.) At $n=622$, indicated by the red arrow, the value of $\Pi^*(n) - e_{max}(n))$ surpasses 0 and never falls below 0 again.}
  \label{fig:Graphs3}
\end{figure}

To verify that no mistakes were made, a graph of  $\Pi(n)$ and  $\Pi^*(n)$ vs. $n$ was made for $n$ not divisible by a prime number less than $\sqrt{n}$. As can be seen, the values of  $\Pi(n)$ lie comfortably within  $\Pi^*(n) \pm e_{max}(n)$ (Figure \ref{fig:Graphs3}). Also, the value of $\Pi^*(n) - e_{max}(n)$ is greater than 0 for all $n \geq 622$. Since $e_{max}(n)$ increases at a much slower rate than $\Pi^*(n)$, $\Pi^*(n) - e_{max}(n)$ will never go below 0 beyond this point. Since it has been shown that Goldbach's conjucture holds true for $n<=622$, this proves that Goldbach's conjecture holds true for all values of even integers $n$.

\section{Summary}

In this paper, it was shown that even integers $n$ that are not divisible by prime numbers less than the $\sqrt{n}$ represent the lower bound on the number of prime pairs of  $n$. The following equation was derived that approximates the number of prime numbers of $n$ for values of $n$ that are not divisible by a prime number less that $\sqrt{n}$.
\begin{flalign*}
\Pi^*(n) = |P(n)|W(n)
\end{flalign*}
where $|P(n)|$ is the number of pairs $(x,y)$ such that $x+y=n$ and $W(x)$ is an approximation to the fraction of prime pairs of $n$ and is defined as follows:
\begin{flalign*}
W(n) = \prod_{p=3}^{l(\sqrt{n})} \frac{(p-2)}{p}
\end{flalign*}
where $l(x)$ is the largest prime number less than $x$ and the product is over prime numbers.\\
It was shown by proof by induction, that  $\Pi^*(n)$ is always greater than 1 and $\Pi^*(n)$ increases indefinitely as $n$ increases.
It was also shown that the maximum error, $e_{max}(n)$, between $\Pi^*(n)$ and the actual number of prime pairs of $n$, $\Pi(n)$, is relatively small and that all values of $n > 622$, $\Pi^*(n) - e_{max}(n) > 0$.

\section{Future Directions}
Future work will involve applying this technique of pairing numbers to prove the Twin Prime Conjecture and Polignac's Conjecture \cite{Polignac}. Polignac's Conjecture states that there is an infinite number of prime pairs $(p_{i},p_{i+1})$ such that $|p_{i+1} - p_{i}| = 2k$ where $k$ is an integer greater than 0. The Twin Prime Conjecture is the case where $k  = 1$.\\
To prove the Twin Prime conjecture, an approximation to the number of twin primes less than $n$, $\pi_{2}(n)$, must be derived.  To do this, odd numbers are paired $(x,y)$ such that $x+2=y$ and $y <= n$. For example, 
(3,5),(5,7),(7,9),(9,11)...,(n-4,n-2),(n-2,n).
Then by eliminating pairs that are divisible by 3, 5, 7, 11 etc, the remaining pairs are twin primes.\\
The number of twin primes less than $n$ will approach the following equation as $n$ gets large:
\begin{center}
$\pi_{2}(n) = |P(n)|W(n)$
\end{center}
This equation is similar to equation \ref{eq:pistar} of the Goldbach proof. This demonstrates that for large values of $n$, the number of twin primes less than $n$ will approach the number of prime pairs of $n$ demonstrating the close relationship between the Twin Prime Conjecture and Goldbach's Conjecture.\\
For other cases of Polignac's Conjecture, for example primes separated by 6, 10 or 30, are cases of the Goldbach Conjecture for $n=6p, n=10p$ or $n=30p$. Thus, Polignac's Conjecture can be proven.\\
Applying this technique to other prime number conjectures will lead to further proofs.

\section{Copyright Notice}
This document is protected by U.S. and International copyright laws. Reproduction and distribution of this document or any part thereof without written permission by the author (Kenneth A. Watanabe) is strictly prohibited.\\
-  I (Kenneth A. Watanabe) grant arXiv.org a perpetual, non-exclusive license to distribute this article.\\ 
-  I (Kenneth A. Watanabe) certify that I have the right to grant this license. \\
-  I (Kenneth A. Watanabe) understand that submissions cannot be completely removed once accepted.\\ 
-  I (Kenneth A. Watanabe) understand that arXiv.org reserves the right to reclassify or reject any submission.\\
Copyright © 2025 by Kenneth A. Watanabe

\end{document}